\def\timestamp{%
Time-stamp: <conjugacyclasses.tex: Friday 11-06-2021 at 11:00:11 (cest)>}
\def\stripname Time-stamp: <#1 #2>{#2}
\edef\filedate{\expandafter\stripname\timestamp}

\documentclass[a4paper]{amsart}

\DeclareMathSymbol\Gap0{AMSb}{`G}
\DeclareMathSymbol\HH0{AMSb}{`H}
\DeclareMathSymbol\I0{AMSb}{`I}
\DeclareMathSymbol\N0{AMSb}{`N}
\DeclareMathSymbol\Q0{AMSb}{`Q}
\DeclareMathSymbol\R0{AMSb}{`R}
\DeclareMathSymbol\Z0{AMSb}{`Z}
\newcommand\betaN{\beta\N}
\newcommand\Nstar{\N^*}

\DeclareMathSymbol\restr\mathbin{AMSa}{"16}
\DeclareMathSymbol\le3{AMSa}{"36}
\DeclareMathSymbol\ge3{AMSa}{"3E}

\newcommand\omegaseq[1]{\langle{#1}_n:n<\omega\rangle}

\newcommand\Triv{\mathsf{Triv}}
\newcommand\Aut{\mathsf{Aut}}

\newcommand\dom{\operatorname{dom}}
\newcommand\ran{\operatorname{ran}}

\newcommand\fin{\mathrm{fin}}

\newcommand\axiom{\mathsf}
\newcommand\CH{\axiom{CH}}

\newcommand\ZFC{\axiom{ZFC}}

\newcommand\cee{\mathfrak{c}}

\newcommand\cl{\operatorname{cl}}
\newcommand\card[1]{\mathopen|{#1}\mathclose|}
\newcommand\orpr[2]{\langle{#1},{#2}\rangle}

\newcommand\pri[1]{}

\usepackage{amsrefs}

\theoremstyle{plain}
\newtheorem{theorem}{Theorem}
\theoremstyle{definition}
\newtheorem{question}{Question}

\begin{document}

\title{Conjugacy classes of autohomeomorphisms of $\Nstar$}

\author{Klaas Pieter Hart}

\address{Faculty EEMCS\\TU Delft\\
         Postbus 5031\\2600~GA {} Delft\\the Netherlands}
\email{k.p.hart@tudelft.nl}
\urladdr{http://fa.ewi.tudelft.nl/\~{}hart}

\author{Jan van Mill}
\address{KdV Institute for Mathematics\\
         University of Amsterdam\\
         P.O. Box 94248\\
         1090~GE {} Amsterdam\\
         The Netherlands}
\email{j.vanmill@uva.nl}
\date{\filedate}

\begin{abstract}
We present some problems related to the conjugacy classes of $\Aut(\Nstar)$.  
\end{abstract}

\maketitle

\section*{Introduction}

\section{Some definitions and notation}

As this note is about the autohomeomorphisms of~$\Nstar$ we fix some notation
regarding~$\betaN$.
For a quick overview of~$\betaN$ we refer to Chapter~D-18 of~\cite{MR2049453};
a more comprehensive introduction is~\cite{MR776630} by the second
author.

We let $\Aut$ denote the autohomeomorphism group of~$\Nstar$,
rather than $\Aut(\Nstar)$, because $\Nstar$~will be the only space under
discussion in this paper.

Let us first identify some easily described members of this group.

\subsection*{Trivial autohomeomorphisms}

To begin: it is clear that an autohomeomorphism of~$\betaN$ leaves 
both~$\N$ and~$\Nstar$
invariant and hence is determined by its restriction on~$\N$, which
is a permutation of~$\N$.
This provides us with the first source of autohomeomorphisms of~$\Nstar$:
the permutation group~$S_\N$ of~$\N$.

For $\pi\in S_\N$ we let $\beta\pi$ denote its extension to~$\betaN$
and $\pi^*$ the restriction of~$\beta\pi$ to~$\Nstar$. 
Thus, permutations of~$\N$ determine autohomeomorphisms of~$\Nstar$.
It is an elementary exercise to show that $\pi^*=\rho^*$ if and only
if the set $\{n:\pi(n)\neq\rho(n)\}$ is finite.
This identifies our first set of easily described members of~$\Aut$:
the image $\{\pi^*:\pi\in S_\N\}$ under the homomorphism $\pi\mapsto\pi^*$.

Every permutation is built up from cyclic permutations and if two permutations,
$\sigma$ and $\tau$ are conjugate, say $\sigma=\pi^{-1}\tau\pi$ then
the permutation~$\pi$ provides a one-to-one correspondents between
the sets of cycles of~$\sigma$ and~$\tau$.
Note that there may also be infinite cycles; 
these look like the infinite cyclic group~$\Z$ with the map~$n\mapsto n+1$.

This shows that conjugacy classes in~$S_\N$ are determined by sequences of
the form $\omegaseq\kappa$, where $\kappa_n$ is the number of~$n$-cycles
in the permutation if $n\ge1$, and $\kappa_0$~is the number of infinite
cycles.
Of course $\kappa_n\le\aleph_0$ for all~$n$.  

There are other autohomeomorphisms with an easy description.
Every bijection $\varphi:A\to B$ between co-finite subsets of~$\N$
determines an autohomeomorphism of~$\N^*$: the restriction~$\varphi^*$
of $\beta\varphi:\cl A\to\cl B$ is a homeomorphism from~$A^*=\N^*$ 
to~$B^*=\N^*$.
As above, if $\psi:C\to D$ is another such bijection then $\varphi^*=\psi^*$
iff $\{n\in A\cap C:\varphi(n)=\psi(n)\}$ is co-finite in~$\N$.

The autohomeomorphisms that we described thus far are called 
\emph{trivial autohomeomorphisms}, they form a subgroup of~$\Aut$
that we will denote~$\Triv$. 
 
Shelah proved that it is consistent that all autohomeomorphisms
of~$\Nstar$ are trivial, see~\cite{MR1623206}*{IV\,\S5}.

\section{Moderately easy results}
\label{sec.easy}

In this section we describe two situations where one can say quite a lot
about conjugacy classes in~$\Aut$.
These are at the opposite ends of the spectrum:
one is the situation where all autohomeomorphisms are trivial
and the other is where the Continuum Hypothesis holds and there is
a wide (possibly the widest) collection of non-trivial autohomeomorphisms
of~$\Nstar$.

\subsection{Trivial autohomeomorphisms}

As mentioned above it is consistent that all autohomeomorphisms are
trivial hence we should look at conjugacy in this case.

We start by quoting a result by Van Douwen from~\cite{MR1035463}.
To this end we associate an integer with every element of~$\Triv$.
Let $\varphi$ be a bijection between co-finite subsets of~$\N$;
define
$$
h(\varphi) = \card{\N\setminus\ran\varphi}-
              \card{\N\setminus\dom\varphi}
$$ 
Now Theorem~6.1 from~\cite{MR1035463} states that $h$~induces a homomorphism
from~$\Triv$ onto~$\Z$.
That is, if $\varphi^*=\psi^*$ then $h(\varphi)=h(\psi)$ and the induced
map $\varphi^*\mapsto h(\varphi)$ is a homomorphism.
We use~$h$ to denote this homomorphism.

If $\varphi^*$ and $\psi^*$ are conjugate \emph{in $\Triv$} 
then $h(\varphi^*)=h(\psi^*)$.
Therefore we concentrate on conjugacy of autohomeomorphisms determined
by members of~$S_\N$.

\subsubsection*{Many conjugacy classes}

Assume $\varphi^*$ and $\psi^*$ are conjugate in $\Triv$, this means that
there is a bijection $\tau:A\to B$ between co-finite sets such that
$\varphi^*\tau^*=\tau^*\psi^*$ and this in turn means that the set
$X=\bigl\{n:\varphi(\tau(n))=\tau(\psi(n))\bigr\}$ is co-finite.

There are only finitely many cycles in $\varphi$ and $\psi$ whose domains
(and their (pre)images under~$\tau$) meet the complement of~$X$.
The remaining cycles of~$\varphi$ are maped by $\tau$ to cycles of~$\psi$
and vice versa. 
Therefore the sequences $\omegaseq\kappa$ and $\omegaseq\lambda$ of 
cycle numbers of $\varphi$ and $\psi$ respectively are almost equal.

This makes it easy to construct a family of $\cee$ many permutations that
represent members of~$\Triv$ that are not conjugate.
For every infinite subset~$x$ of~$\N$ take a partition $\{A_n:n\in x\}$ 
of~$\N$ such that $\card{A_n}=n$ for all~$n$ and create a permutation~$\pi_x$ 
of~$\N$ by permuting each~$A_n$ cyclically ---say $(a_1\,a_2\,\ldots\,a_n)$, 
where $A_n=\{a_1,a_2,\ldots,a_n\}$ listed in order.
 
If $x\neq y$ then $\pi_x$ and $\pi_y$ are not conjugate in~$S_\N$,
but they may of course be conjugate in~$\Triv$, say if
$x=\{2,5\}\cup\{n:n\ge10\}$ and 
$y=\{3,4\}\cup\{n:n\ge10\}$.
If the symmetric difference of $x$ and $y$ is infinite then $\pi_x^*$ 
and~$\pi_y^*$ will not be conjugate.
This implies that an almost disjoint family of cardinality~$\cee$ will provide
us with $\cee$~many conjugacy classes.

In Section~\ref{sec.questions} we raise some questions suggested by these
considerations.

\subsubsection*{Infinite cycles}

We should make a few remarks about infinite cycles in permutations
of~$\N$.
Such a cycle is, as mentioned above, a copy of the set~$\Z$ is integers
with the shift map $\sigma:n\mapsto n+1$.
If we work, for the moment, in $\beta\Z$ then we see that $\Z^*$
is split into two clopen sets that are minimally $\sigma^*$-invariant.
Indeed, it should be clear that $L=\{n\in\Z:n<0\}^*$ and $R=\{n\in\Z:n\ge0\}^*$
are both invariant under~$\sigma^*$.
It is only slightly more difficult to verify that if $A$ is an infinite
subset of~$L$ (or~$R$) such that $\sigma*[A^*]\subseteq A^*$
then $L\setminus A$ (or $R\setminus A$) is finite.

We see that when passing from $S_\N$ to $\Triv$ an infinite cycle ceases
to be a unit: it splits into two independent autohomeomorphisms.

Any conjugation, even if non-trivial, will preserve the structure
of these minimal invariant clopen subsets; in particular the cardinality
of the family of these sets.
This shows that for two permutations $\varphi$ and $\psi$ for which
$\varphi^*$ and $\psi^*$ are conjugate in~$\Aut$ the numbers of infinite
cycles are be the same. 

Conclusion: if $\varphi$ and $\psi$ have cycle number sequences
$\omegaseq\kappa$ and $\omegaseq\lambda$ respectively and if $\varphi^*$
and $\psi^*$ are conjugate in~$\Triv$ then 
\begin{itemize}
\item $\kappa_0=\lambda_0$,
\item for all $n\ge1$ the equalities $\kappa_n=\aleph_0$ and 
      $\lambda_n=\aleph_0$ are equivalent, and  
\item $\kappa_n=\lambda_n$ for all but finitely many~$n\ge1$.
\end{itemize}
The last condition becomes important only if there are infinitely
many~$n$ for which $\kappa_n$ and $\lambda_n$ are finite.

Part of the analysis above was used in~\cite{MR3563083} to show that 
an autohomeomorphism of~$\Nstar$ derived from a homeomorphism between~$\Nstar$
and~$\omega_1^*$ was non-trivial.

\subsection{The Continuum Hypothesis}

Many questions have a relatively easy answer under the assumption of~$\CH$.
This is largely due to Parovichenko's characterization of~$\Nstar$
under that assumption.

The Continuum Hypothesis implies that $\Aut$~is a simple group.
This was proven by Fuchino in~\cite{fuchino-thesis} in a more
general form: the automorphism group of a saturated Boolean algebra
is simple, see~\cite{MR991606}*{Theorem~5.12} for a more accessible proof. 
Since the Continuum Hypothesis implies that the Boolean algebra
of clopen sets of~$\Nstar$ is saturated the result follows.
In~\cite{MR1164729} Fuchino proved that $\Aut$~is also simple in the
$\aleph_2$-Cohen model.

We shall show that $\CH$ implies that $\Aut$~has $2^\cee$ many conjugacy 
classes.
For this we need two known results about~$\Nstar$.

The first result is due to Hart and Vermeer.

\begin{theorem}[\cite{MR1260168}, $\CH$]\label{thm.P-set}
Every $P$-set in $\Nstar$ is the fixed-point set of an involution.  \qed
\end{theorem}

In fact the proof is flexible enough to enable one to make any given $P$-set
the fixed-point set of an autohomeomorphism of any prescribed finite order.

\smallbreak
The next result is due to Dow, Gubbi and Szyma\'nski.

\begin{theorem}[\cite{MR929014}]
There are $2^\cee$ many mutually non-homeomorphic (rigid)
separable extremally disconnected spaces.  \qed
\end{theorem}

We combine these two results using the well-known fact that under $\CH$ every
such separable space can be embedded into~$\Nstar$ as a $P$-set,
see~\cite{MR776630}*{Theorem~1.4.4}.

This produces $2^\cee$ many mutually non-homeomorphic $P$-sets.
Each of these is the fixed-point set of an involution.
These involutions are never conjugated because conjugate autohomeomorphisms
have homeomorphic fixed-point sets.

Although this set of involutions answers the question about the number of 
conjugacy classes of~$\Aut$ it is actually quite small.
As noted above we can, almost for free, get autohomeomorphisms of any desired
finite order.
In addition, Theorem~1.4.4 from~\cite{MR776630} states that every compact
$F$-space of weight~$\cee$ can be embedded in~$\Nstar$ as a nowhere 
dense $P$-set.
This immediately gives us many more conjugacy classes.

It also suggests some questions that we shall mention in 
Section~\ref{sec.questions}.

\section{Questions}
\label{sec.questions}

In this section we collect questions that are suggested by the results
in Section~\ref{sec.easy} and by other results in the literature.

\subsubsection*{What happens to trivial autohomeomorphisms?}

We found $\cee$ many conjugacy classes in~$S_\N$ by exploiting
the cycle structure of permutations.
In the model where all members of~$\Aut$ are trivial these gave
us the maximum possible number of conjugacy classes.

We have also seen that permutations that are not conjugate may induce
the same autohomeomorphism of~$\Nstar$.

It also seems conceivable that quite distinct permutations may determine
conjugate members of~$\Aut$ in case there are non-trivial
autohomeomorphisms.
The general question then is: what happend to conjugacy classes of
trivial autohomeomorphisms when $\Aut$ is not equal to~$\Triv$?

We give some specific versions of this question below, where we should 
emphasize that to the best of our knowledge these questions have not even 
been answered under the assumption of the Continuum Hypothesis
when $\Aut$~is much much richer than~$\Triv$.

\begin{question}
What is the relationship between conjugacy classes of permutations in~$S_\N$
and their conjugacy classes in~$\Aut$?  
\end{question}
  
This question is quite general and we may specialize to the permutations
we considered in Section~\ref{sec.easy}.

\begin{question}
Let $x$ and $y$ be infinite subsets of~$\N$ such that $x\neq^*y$.
Under what conditions will $\pi_x^*$ and $\pi_y^*$ become conjugate
in~$\Aut$?  
\end{question}

As toy problems one may consider $x=\{2^n:n\in\N\}$ and $y=\{3^n:n\in\N\}$,
or $u=\{2^n:n\in\N\}$ and $v=\{4^n:n\in\N\}$.
 
As a variation we can use a function $f:\N\to\N$ to specify a 
permutation~$\tau_f$ (up to conjugacy):
partition $\N$ into set $A_n$ where $\card{A_n}=f(n)$ for all~$n$ and
turn each~$A_n$ into an $f(n)$-cycle.
The difference with the $\pi_x$ is that we allow repetitions of cardinalities.

Once one knows the effect of~$\CH$ on these questions one can venture
into models where there are autohomeomorphisms of varying degrees
of (non-)triviality.
A sample of such models can be found for example 
in~\cites{MR2437015,MR2506596,MR1002627,MR1271551}

\subsubsection*{Infinite cycles}

In the above questions we concentrated on finite cycles.
We have seen that an infinite cycle ceases to be a building block
when we move to $\Aut$.
It gives us two autohomeomorphisms that are not induced by permutations:
\begin{itemize}
\item the forward shift $\sigma_N:n\mapsto n+1$ on $\N$, and
\item the downward shift $\sigma_\N^{-1}: n\mapsto n-1$ on $\N$.
\end{itemize}
The latter two are mapped to $1$ and $-1$ respectively by the homomorphism~$h$.

The two shifts are minimal in that $\Nstar$ (and the empty set) are the only
clopen sets that are invariant.
The two shifts are not conjugate in~$\Triv$, but whether they can be
conjugate is open, even under~$\CH$.

\begin{question}
Is it consistent that $\sigma_\N$ and $\sigma_N^{-1}$ are conjugate?  
\end{question}

An extensive study of this problem can be found in~\cite{MR4138425}.

The shift map also has various universality properties, 
see~\cites{MR3835074,MR4013982}; for example $\CH$~implies that the
system $\orpr\Nstar{\sigma_\N^{-1}}$ is a quotient of~$\orpr\Nstar{\sigma_\N}$.

\subsubsection*{Other ways of (dis)proving conjugacy}

We exhibited, under~$\CH$, many conjugacy classes by exhibiting 
autohomeomorphisms with non-homeomorphic fixed-point sets.
These fixed-point sets were all $P$-sets and that is no coincidence;
the converse of Theorem~\ref{thm.P-set} is a theorem of~$\ZFC$: 
every fixed-point set of an autohomeomorphism of~$\Nstar$ is a $P$-set.

This, combined with the homeomorphism extension theorem for $P$-sets
from~\cite{MR1277871}, indicates that the fixed-points sets will play a key
role in deciding conjugacy. 

As we saw above a $P$-set can be the fixed-point set of autohomeomorphisms 
of all possible finite orders.

This suggests the following concrete problem, under $\CH$:

\begin{question}
Assume $h$ and $g$ are two autohomeomorphisms with the same fixed-point set
and the same finite order.
Are $h$ and $g$ conjugate?  
\end{question}

We end with a general question: what other invariants can we use
to (dis)prove conjugacy of autohomeomorphisms of~$\Nstar$.

\end{document}